\documentclass[12pt]{amsart}
\usepackage{amssymb,hyperref,graphicx,epsfig,color,bm,amsmath} 
\usepackage{float}

\newtheorem{theorem}{Theorem}

\newtheorem{conjecture}[theorem]{Conjecture}
\newtheorem{corollary}[theorem]{Corollary}
\newtheorem{definition}[theorem]{Definition}
\newtheorem{example}[theorem]{Example}

\newtheorem{exercise}[theorem]{Exercise}
\newtheorem{fact}[theorem]{Fact}
\newtheorem{lemma}[theorem]{Lemma}
\newtheorem{problem}[theorem]{Problem}
\newtheorem{proposition}[theorem]{Proposition}
\newtheorem{question}[theorem]{Question}
\newtheorem{remark}[theorem]{Remark}

\newcommand{\bcon}{\begin{conjecture}}
\newcommand{\econ}{\end{conjecture}}
\newcommand{\bcor}{\begin{corollary}}
\newcommand{\ecor}{\end{corollary}}
\newcommand{\bdf}{\begin{definition}}
\newcommand{\edf}{\end{definition}}
\newcommand{\beq}{\begin{equation}}
\newcommand{\eeq}{\end{equation}}
\newcommand{\bexa}{\begin{example}}
\newcommand{\eexa}{\end{example}}
\newcommand{\bexe}{\begin{exercise}}
\newcommand{\eexe}{\end{exercise}}
\newcommand{\bfac}{\begin{fact}}
\newcommand{\efac}{\end{fact}}
\newcommand{\bite}{\begin{itemize}}
\newcommand{\eite}{\end{itemize}}
\newcommand{\blem}{\begin{lemma}}
\newcommand{\elem}{\end{lemma}}
\newcommand{\bprb}{\begin{problem}}
\newcommand{\eprb}{\end{problem}}
\newcommand{\bpro}{\begin{proposition}}
\newcommand{\epro}{\end{proposition}}

\newcommand{\bque}{\begin{question}}
\newcommand{\eque}{\end{question}}
\newcommand{\brem}{\begin{remark}}
\newcommand{\erem}{\end{remark}}
\newcommand{\bthm}{\begin{theorem}}
\newcommand{\ethm}{\end{theorem}}
\newcommand{\bmat}{\begin{matrix}}
\newcommand{\emat}{\end{matrix}}

\newcommand{\bpr}{\begin{proof}}
\newcommand{\epr}{\end{proof}}

\newcommand{\lb}{\label}

\newcommand{\com}[1]{\,}

\newcommand{\cal}{\mathcal}
\newcommand{\p}{\partial}

\newcommand{\Z}{\mathbb Z}

\newcommand{\C}{\mathbb C}

\newcommand{\ve}{\varepsilon}

\newcommand{\diag}[2]{\parbox{#2}{\includegraphics[width=#2]{#1.pdf}}} 

\title{SU(3)-skein algebras and webs on surfaces}

\author{Charles Frohman, Adam S. Sikora}


\begin{document}

\thispagestyle{empty}

\begin{abstract} 
The $SU_3$-skein algebra of a surface $F$ is  spanned by isotopy classes of  certain framed graphs in $F\times I$ called $3$-webs subject to the skein relations encapsulating relations between $U_q(sl(3))$-representations. 
It is expected that their theory parallels that of the Kauffman bracket skein algebras. We make the first step towards developing that theory by proving that the reduced $SU_3$-skein algebra of any surface of finite type is finitely generated.

We achieve that result by developing a theory of canonical forms of webs in surfaces. Specifically, we show that for any ideal triangulation of $F$ every reduced $3$-web 
can be uniquely decomposed into unions of pyramid formations of hexagons and disjoint arcs in the faces of the triangulation with possible  
additional ``crossbars'' connecting their edges along the ideal triangulation. We show that such canonical position is unique up to ``crossbar moves''.
That leads us to an associated system of coordinates for webs in triangulated surfaces 
(counting intersections of the web with the edges of the triangulation and their rotation numbers inside of the faces of the triangulation) which determine a reduced web uniquely. 

Finally, we relate our skein algebras to $\cal A$-varieties of Fock-Goncharov and to $\text{Loc}_{SL(3)}$-varieties of Goncharov-Shen. We believe that our coordinate system for webs is a manifestation of a (quantum) mirror symmetry conjectured by Goncharov-Shen.
\end{abstract}

\address{14 MacLean Hall, University of Iowa, Iowa City, Iowa 52242, and 244 Math Bldg, University at Buffalo, SUNY, Buffalo, NY 14260}
\email{charles-frohman@uiowa.edu, asikora@buffalo.edu}

\pagestyle{myheadings}

\maketitle

\tableofcontents{}


%
\section{Introduction}
%

The $SU(n)$-skein modules of $3$-manifolds $M$ were introduced in \cite{S2}. They are built of $n$-valent graphs, called $n$-webs, considered up to skein relations encapsulating relations appearing in the $U_q(sl(n))$-representation theory. By \cite{S2}, they are $q$-deformations of the $SL(n)$-character variety of $\pi_1(M)$,  \cite{S2}.

For the cylinder over a surface, $F\times I,$ the skein module forms an algebra denoted by $S_n(F).$ 
For $n=2$ this is the Kauffman bracket skein algebra, which is one of the central objects of quantum topology, and has been a subject of a very active research over the recent years -- see for example \cite{BW1, BW2, BFK, FrGe, FKL, Ma, Le1, Mu, PS1, PS2, Th} and references within. That research connected skein algebras to quantum invariants of $3$-manifolds, (quantum) hyperbolic geometry, (quantum) cluster algebras, (derived algebraic geometry) of moduli spaces, and the topological quantum field theory.

In this paper, we make the first step towards a development of an analogous theory for $SU(3)$-skein algebras, building on \cite{SW} and \cite{Ku}. 

Since the combinatorics associated with isotopy of webs in surfaces is much richer than for embedded one-manifolds, $SU(n)$-skein algebras are much harder to analyze for $n>2$ than for $n=2.$ 

From now we will work with $SU(3)$-skein algebras only, refer to $3$-webs as ``webs'', and allow for web endpoints in $\p F$. To be specific we will be working mostly with the reduced version, $\cal{RS}(F,M),$ of this algebra obtained by killing all webs with boundary parallel 1-, 2- and 3-gons.

The main results of this paper are:
\bite
\item Theorem \ref{t-coordinates} which for marked surfaces $(F,M)$ with an ideal triangulation provides an explicit basis of $\cal{RS}(F,M)$ of ``canonical webs" and provides a coordinate system for those webs. 
\item  Theorem \ref{t-main} stating that the reduced $SU(3)$-skein algebras of all marked surfaces of finite type are finitely generated. (An analogous result for the Kauffman bracket skein relations without marked points is due to Bullock, \cite{Bu} (cf. also \cite{AF}) and with marked points in \cite{PS2}.)  
\eite

Let us elaborate on the first result. We show that each reduced web (i.e. a non-elliptic web without boundary parallel 1-, 2- and 3-gons) admits a canonical position through ``tidying up'', Theorem \ref{t-main-canonical}, in which it is composed of a disjoint union of certain hexagonal patterns (which we call pyramids) and some simple arcs in each face of the triangulation. These triangle pieces may have additional ``crossbars'' connecting their edges along the edges of ideal triangulation. We show that such canonical position is unique up to ``crossbar moves''.
That leads us to an associated system of coordinates for webs in triangulated surfaces 
(counting intersections of the web with the edges of the triangulation and their rotation numbers inside of the faces of the triangulation) which determine a reduced web uniquely (Theorem \ref{t-coordinates}). An alternative but similar system of coordinates can be found in recent \cite{DS1, DS2}. (This is an updated version of our original manuscript posted on arXiv in February 2020, with an expanded discussion of the connections of our work to Fock-Goncharov theory. In particular, we reference here some of the preprints posted since then.)

We consider this work as part of higher quantum Teichm\"uller theory. In particular, we believe that our reduced skein algebra is a quantization of $Loc_{SL(3,\C),F}$-variety of  \cite{GS1}, cf. Conjecture \ref{c-phi3}. (For $F$ with no unmarked boundary components, $Loc_{SL(3,\C),F}$ coincides with the $\cal A$-variety of \cite{FoGo}.)

Assuming the above conjecture our Proposition \ref{pro-basis-RS} defines a basis of 
$\cal{RS}(F,M)=\C[Loc_{SL(3,\C),F}]$ equivariant under the mapping class group, whose existence was conjectured by Goncharov-Shen, \cite{GS1}.
Their mirror-duality conjecture asserts that this basis is in bijection with positive, tropical points of the Fock-Goncharov $\cal A$-variety, $\cal A_{PSL(3,\C),F}^+(\Z^t)$. We believe that our coordinate system on reduced non-elliptic webs is a manifestation of that duality. See further discussion in Sec. \ref{s-FG}.



%
\section{Acknowledgements}
%

We would like to thank Daniel Douglas, Zhe Sun and the referee of this paper for their insightful comments about the relation of this work to Fock-Goncharov theory.

%
\section{$SU_3$-skein algebras}
\lb{s-skein-alg}
%

A \underline{marked surface} is an oriented surface $F$ together with a finite subset $M$ of $\p F$ of \underline{marked points}. ($\p F$ may be empty). In Sec. \ref{s-ideal}, we will alternatively consider $\hat F=F-M$ instead with a hyperbolic metric on it so that $M$ is a set of points at infinity and the boundary arcs of $F$ connecting them are infinite geodesics. 

Let $I=[0,1].$ A \underline{web in $(F,M)\times I$}  consists of
\bite
\item $1$-valent external vertices in fibers $\{m\}\times I$ over (some of) the marked points $m\in M$
\item $3$-valent internal vertices which are either sources or sinks.
\item oriented edges connecting the (internal and external) vertices
\item oriented loops
\item A framing, which can be formally defined by
\bite
\item  a choice of a properly  embedded oriented   compact surface $S$ having the web as a spine, called framing, such that the projection on the first factor $p:F\times I\rightarrow F$ restricts to a local homeomorphism 
$p|_S:S\rightarrow F$.
\item $S\cap \partial F \times I$ is a collection of horizontal arcs, 
each of which containing a unique point of $M$.
\eite
\eite 

A web may be empty, $\emptyset$, or disconnected. If it does not contain any internal vertices then it is a union of ribbons and annuli. Webs in $(F,M)\times I$ are considered up to \underline{isotopy} which is a homotopy within the space of all such webs.

We say that a boundary component of $F$ is \underline{marked} if it contains a point of $M$.
Connected components of marked boundary components of $F$ with points of $M$ removed are the \underline{boundary intervals} of $(F,M)$. We orient them arbitrarily (for example, with the orientation induced from $F$) and we consider that choice of orientation as a part of the structure of a marked surface $(F,M).$

It is convenient to represent webs by their $1$-dimensional diagrams in $F$. 
A \underline{web diagram} in $(F,M)$ is an oriented graph composed of
\bite
\item $1$-valent (external) vertices in boundary intervals of $(F,M)$ (i.e. away from the marked points)
\item $3$-valent (internal) vertices, sinks or sources, in the interior of $F$
\item crossings, i.e. $4$-valent vertices with the overpass marked, cf. Figure \ref{f-internal-v}.
\item oriented edges connecting the above vertices
\item oriented loops.
\eite

\begin{figure}[h]
\begin{center}
\diag{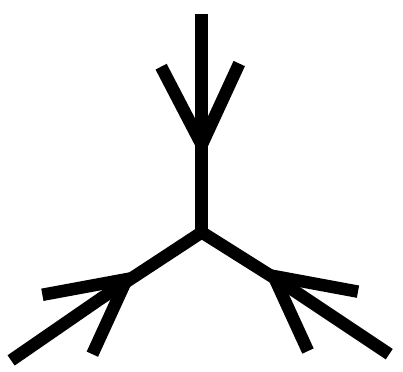}{.4in}\hspace*{.4in} \diag{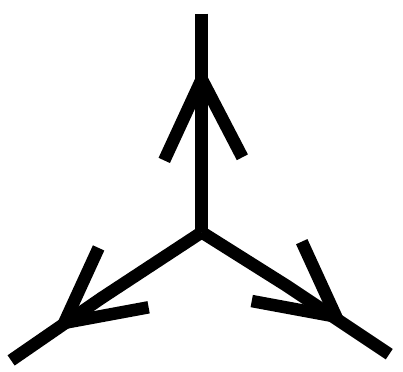}{.4in}\hspace*{.4in}
\diag{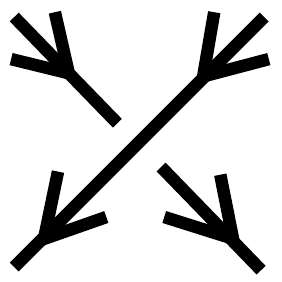}{.3in}\hspace*{.4in} \diag{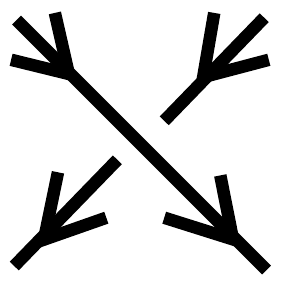}{.3in}
\end{center}
\caption{Types of internal vertices in $3$-web diagrams}
\lb{f-internal-v}
\end{figure}

Each web diagram represents a web in $(M,F)\times I$ obtained by pushing the diagram endpoints (external vertices) in every marked boundary interval starting with $m\in M$ into $\{m\}\times I \subset F\times I$ so that the points closer to $m$ go to lower levels. (Since marked boundary components are oriented, each marked boundary interval has a starting point.) This operation in the context of framed links was discussed in \cite{PS2}. Keeping the web endpoints away from $M$ rather than at $M$ makes the theory easier.

Conversely, every web in $(M,F)\times I$ is represented by a diagram,
which is unique up to the Reidemeister-type moves of Fig. \ref{fig-Reidemeister}.

\begin{figure}[h]
\centerline{\diag{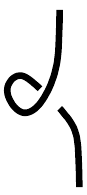}{.2in} $\leftrightarrow$ \diag{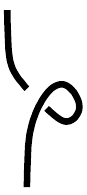}{.2in},\quad
\diag{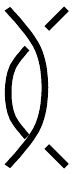}{.2in} $\leftrightarrow$ \diag{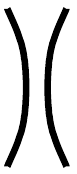}{.2in},\quad
\diag{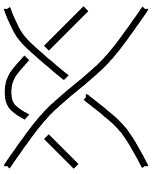}{.4in} $\leftrightarrow$ \diag{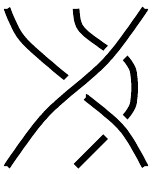}{.4in},\quad
\diag{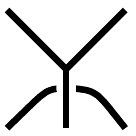}{.4in} $\leftrightarrow$ \diag{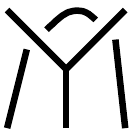}{.4in}}
\caption{Reidemeister-type moves for $3$-web diagrams. Orientations are not marked and all are possible.}
\lb{fig-Reidemeister}
\end{figure}

One of the results of this paper is a system of coordinates for webs in $(F,M)$, cf. Sec. \ref{s-coordinates}.

Let $R$ be a commutative ring with a specified invertible element $q^{1/3}$. 
We use standard exponential notation to denote the powers of this element, eg. $q$ is its cube.
The \underline{$SU(3)$-skein algebra} $\cal S(F,M)$ of $(F,M)$ is the free $R$-module with the basis of all webs in $(F,M)\times I$
quotiented by the following skein relations: 
\beq\lb{skein1} \diag{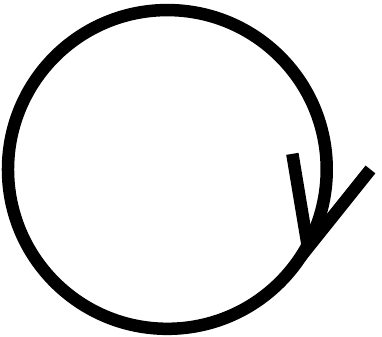}{.3in}\ -\ (q^2+1+q^{-2})\emptyset,\quad
\diag{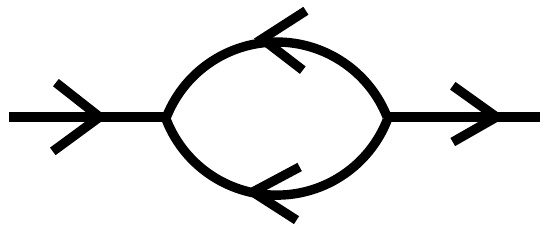}{.7in}+
    (q+q^{-1})
\diag{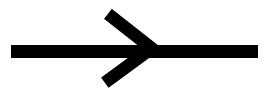}{.3in}
\eeq
\beq\lb{skein2} \diag{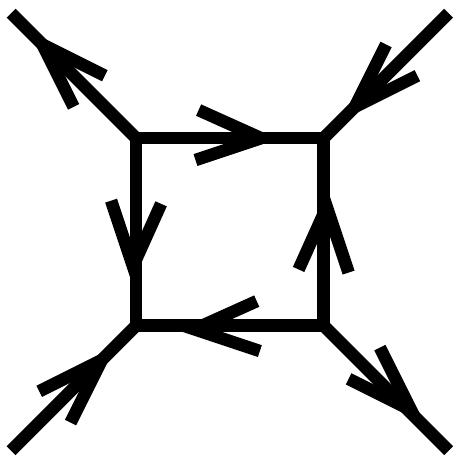}{.45in}-\diag{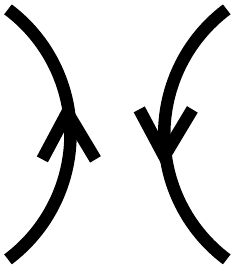}{.3in}-\diag{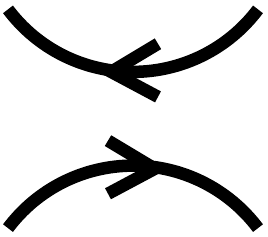}{.35in}
\eeq
\beq\lb{skein3} 
 \diag{crossp}{.3in}-  q^{-\frac{1}{3}}\diag{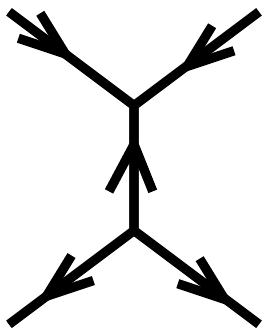}{.35in}-
q^{\frac{2}{3}}\diag{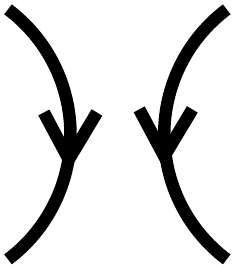}{.3in},\quad \diag{crossm}{.3in} -  q^{\frac{1}{3}}
\diag{crossh}{.35in}-q^{-\frac{2}{3}}\diag{smootho}{.3in}.
\eeq
All these diagrams have the standard, horizontal framing. (Our $q$ is $q^{1/2}$ in \cite{Ku}.) 

Let $W_1\cdot W_2$ denote a web obtained by stacking the web $W_1$ on top of $W_2$ in 
$F\times I$. That operation extends to a product on $\cal S(F,M)$ making it into an associative $R$-algebra with the identity $\emptyset.$ (That explains the term ``algebra'' in its name.)  The relation of this algebra to the stated skein algebras of \cite{Hi,LS} is discussed in \cite{LS}.

As alluded earlier, it is easier to operate on web diagrams than webs. 
In the diagrammatic approach, the stacking  of $W_1$ on top of  $W_2$ requires that the endpoints of $W_1$ lie after the points of $W_2$ in the respective boundary intervals of $(F,M).$ (Here again we use the fact that boundary intervals are oriented. Kauffman bracket skein algebras for surfaces with marked points were considered in 
\cite{BW1, Le1, Le2, Mu, PS2}. Our approach follows that last paper.)

Note that the above skein relations allow for resolving any crossing as a linear combination of
the \underline{vertexless} and the \underline{I-resolutions}, Fig. \ref{fig-cross-res},

\begin{figure}[h]
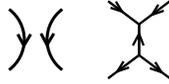

\centerline{\diag{smootho}{.3in} \quad \diag{crossh}{.35in}}
\caption{The vertexless resolution and the I-resolution of a crossing.}
\lb{fig-cross-res}
\end{figure}

\bdf\lb{r-iso-crossingless} 
We say that two web diagrams are \underline{isotopic} if they represent webs isotopic in $(F,M)\times I.$ We say that two crossingless web diagrams are \underline{planarly isotopic} if they are isotopic in the space of crossingless web diagrams in $(F,M).$
\edf

Note that neither isotopy no planar isotopy allows for endpoints passing each other in $\p F.$
Isotopy is more flexible than planar isotopy, because it allows for a flip of two parallel loop components, (made up of two Reidemeister II moves), cf. Fig. \ref{f-isotopic}. We call it a \underline{flip move}.
\begin{figure}[h]
\centerline{\diag{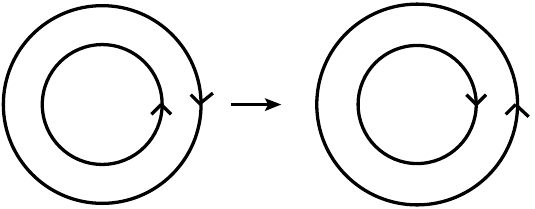}{1.7in}}
\caption{A flip move.}
\lb{f-isotopic}
\end{figure}
In fact, any two isotopic crossingless web diagrams differ by a planar isotopy and flip moves.

A crossingless web diagram without trivial components (contractible loops) nor internal  2- or 4-gons is called \underline{non-elliptic}. (``Internal'' refers to regions in the interior of $F$.)

The results and methods of \cite{SW} imply:

\bthm\lb{thm-basis}
Non-elliptic web diagrams in $(F,M)$ up to isotopy form a basis of $\cal S(F,M).$
\ethm

We have seen in \cite[Example 8]{PS2} that the Kauffman bracket skein algebra is not finitely generated in general, even in the case of $(D^2,M)$ for $M\ne \emptyset.$
For the same reason, the $SU_3$-skein algebra of $(D^2,M)$ is not finitely generated either.
To achieve a finite generation of the $SU(3)$-skein algebras we will consider their reduced version, analogous to that in \cite{Le1,Le2,Mu, PS2}.

%
\section{Reduced $SU_3$-skein algebras}
\lb{s-reduced-skein}
%

Let $a$ be a fixed invertible element of $R$ and let $\cal D(F,M)$ be the $R$-submodule of $\cal{S}(F,M)$ spanned by the boundary skein relations: 
\beq\label{e-bskein1}
 \diag{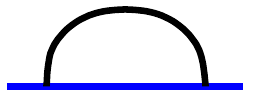}{.6in}= 0= \diag{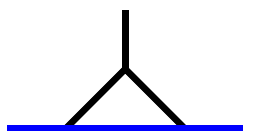}{.6in}, 
\eeq 
\beq\label{e-bskein2}
\diag{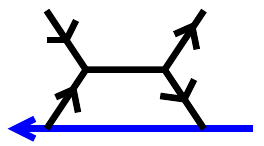}{.6in}= a\cdot \diag{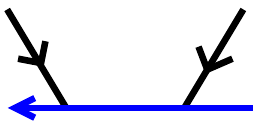}{.6in},\hspace*{.4in} \diag{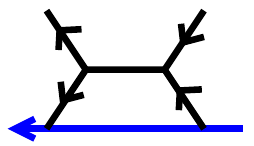}{.6in}= a^{-1}\cdot \diag{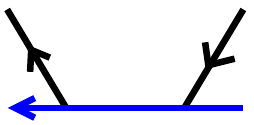}{.6in},
\eeq
where the blue horizontal line denotes $\p F$. The orientations of the arc, the trigon, and the boundary segment 
in \eqref{e-bskein1} are arbitrary. (One can show that each of the relations of \eqref{e-bskein2} is a consequence of the other one and the remaining relations.)

\blem\lb{l-D-ideal}
$\cal D(F,M)$ is a $2$-sided ideal in $\cal S(F,M).$
\elem

\bpr
For any web diagrams $W_1, W_2$ in $(F,M),$ the endpoints of $W_2$ in $W_1\cdot W_2$ and
in $W_2\cdot W_1$ are either after or before the endpoints of $W_1$ in each boundary interval. Consequently, they do no affect the above skein relations.
\epr

$\cal{RS}(F,M)=\cal S(F,M)/\cal D(F,M)$ is the \underline{reduced $SU(3)$-skein algebra} of $(F,M)$ which is the main subject of this paper.


A non-elliptic web diagram is \underline{reduced} if it does not contain any path of $\leq 3$ edges parallel to a boundary interval, cf. \eqref{e-bskein1}-\eqref{e-bskein2}.
(A path of edges of $W$ is parallel to a boundary interval if together with it bounds a $2$-disk
whose interior is disjoint from $W$.)

\begin{figure}[h]
\centerline{\diag{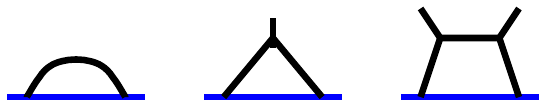}{1.6in}}
\caption{Paths of $1$, $2$ and $3$ edges parallel to a (blue) boundary interval. The orientations are unmarked and arbitrary.}
\lb{fig-degeneracies}
\end{figure}

These web diagrams considered up to isotopy span $\cal{RS}(F,M)$. However, for $M\ne \emptyset$, the relation of Figure \ref{fig-parallel-arcs} shows a linear dependence between them.
\begin{figure}[h]
\centerline{$a^{-k}\cdot$\diag{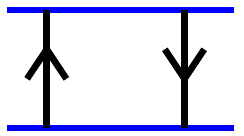}{0.6in}$=$ \diag{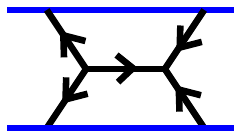}{0.6in}
$=$ $a^k\cdot$ \diag{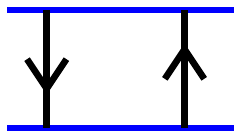}{0.6in}}
\caption{Relations between parallel arcs. $k\in \{-2,0,2\}$ depends on the orientations of the boundary arcs.}
\lb{fig-parallel-arcs}
\end{figure}

For that reason we say that two reduced crossingless web diagrams are \underline{equivalent} if they
are related by isotopy of webs (of Definition \ref{r-iso-crossingless}) and by a permutation of 
parallel arc components.
By the method of confluence, cf. \cite{SW}, we prove:

\bpro\lb{pro-basis-RS}
$\cal{RS}(F,M)$ has a basis consisting of one representative per each equivalence class of reduced crossingless web diagrams in $(F,M)$.
\epro

In order to make such a basis explicit, we need to choose a convention about parallel loops and arcs similar to the conventions governing car traffic.

\begin{figure}[h]
\centerline{\diag{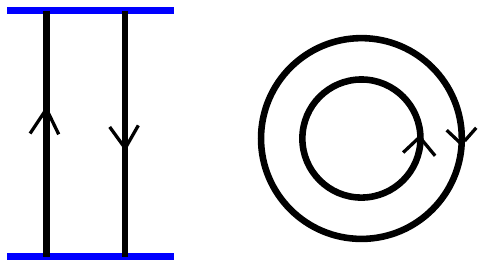}{1.1in}}
\caption{British highways. Blue horizontal arcs belong to $\p F.$}
\lb{fig-British-highways}
\end{figure}

For the sake of precision, all pictures in this paper show $F$ oriented counterclockwise.
We call a pair of parallel loops or arcs a \underline{British highway} if they are oriented inconsistently with the boundary of the annulus or of the rectangle they bound.
Note that such parallel arcs and loops resemble the left-side traffic of the countries of the former British empire. As for the sake of safe and efficient traffic, one needs to mandate either the left-side and the right-side rule, so we do in order to fix a basis of $\cal{RS}(F,M)$.
Specifically, as a consequence Proposition \ref{pro-basis-RS} we obtain:

\bcor\lb{cor-basis-RS}
$\cal{RS}(F,M)$ has a basis of reduced webs in $(F,M)$ (up to planar isotopy) without British highways.
\ecor

Let us discuss generating sets for skein algebras now.

A \underline{triad} is a connected web in $(F,M)\times I$ with a single $3$-valent vertex (source or sink). 

We say that an arc diagram in $(F,M)$ is \underline{descending} if transversing it according to its orientation one passes each crossing through an overpass first. Similarly, a triad is descending if there is an ordering of its three edges such that the above property holds when transversing first the 1st edge, then the 2nd, and finally the 3rd one.
Finally, a knot diagram $D$ is descending if there is a point in $D$ such that transversing $D$ according to its orientation, starting from that point, one passes each crossing through an overpass first.

A knot, arc, or triad in $(F,M)\times I$ is descending if it has a descending diagram.

\bthm[Proof in Sec. \ref{s_t-main-vr-proof}]\lb{t-main-vr}
(1) For every marked surface $(F,M)$, the algebra $\cal R\cal S(F,M)$ is generated by descending knots, arcs, and triads.\\
(2) $\cal R\cal S(D^2,M)$ is generated by crossingless arc and triad diagrams with endpoints at distinct boundary intervals.
\ethm

Since each crossingless arc or triad diagram in $(D^2,M)$ is determined by its orientation and the boundary intervals containing its endpoints, the cardinality of the above generating set for $F=D^2$ is $2{|M|\choose 2}+2{|M|\choose 3}$, where
${k \choose l}=0$ for $k<l.$ However, for $F\ne D^2$ the above generating set is infinite.

The following much stronger and harder to prove result is one of the two main results of this paper. It is a consequence of its graded version, Theorem \ref{main-gr}. Following the standard definition, we say that $(F,M)$ is \underline{of finite type} if $F$ is obtained by removing a finite set of internal points from a compact surface.

\bthm[Proof in Sec. \ref{s_graded}]\lb{t-main}
${\cal RS}(F,M)$ is a finitely generated $R$-algebra for every marked surface $(F,M)$ of finite type.
\ethm

%
\section{Alternative Approach to Marked Surfaces, Ideal triangulations}
\lb{s-ideal}
%

Since ${\cal RS}(F,M)$ is not affected by the removal of unmarked boundary components from $F$, we can assume that they are never there (and think of them as ``punctures'').

Furthermore, by abuse of notation, from now on we will identify $(F,M)$ with $\hat F=F-M$. 

Furthermore, from now on we will not consider webs in $(F,M)\times I$ directly but rather their representations by web diagrams in $(F,M)$ or in $\hat F$ only. (This does not lead to any issues since web diagrams cannot touch points of $M$ anyways.)
For that reason we will call web diagrams \underline{webs in $\hat F$} for simplicity.  

Let us denote the genus of $F$ by $g$ and the number of boundary components and ends of $F$ by $k$. 

\blem\lb{l-hyp} If $4-4g-2k-|M|<0$ then $\hat F$ can be given a hyperbolic metric such that  
\bite
\item points $M$ are at infinity and all boundary intervals are infinite geodesics. 
\item all ends of $\hat F$, other than points of $M$ are punctures.
\eite 
\elem

\bpr  Let $\tilde F$ be two copies of $\hat F$ glued along their corresponding boundary intervals. The Euler characteristic of $\tilde F$ is $4-4g-2k-|M|<0$ and, hence, it admits a complete hyperbolic structure. The boundary intervals of $F$ are arcs in $\tilde F$ connecting its ends and so they can be isotoped to geodesics.
\epr

Although $\hat F$ is not necessarily hyperbolic, we will call its ends other than the points of $M$ \underline{punctures}. An \underline{infinite arc} in $\hat F$ is an embedding of $(-\infty,\infty)$ into $\hat F$ with each of its ends either in $M$ or one of the punctures of $\hat F.$ (In particular, every boundary arc of $(F,M)$ is an infinite arc in $\hat F$.)

An \underline{ideal triangulation} of $\hat F$ is a disjoint collection of infinite arcs $\gamma_1,...,\gamma_N$ in $\hat F$ which includes all boundary intervals of $\hat F$
and which decomposes $\hat F$ into ideal triangles.

By Lemma \ref{l-hyp}:

\bcor\lb{c-triangulable} For every non-closed marked surface $(F,M)$, other than the open disk, 
the closed disk with at most two marked points, or an annulus $(J\times S^{1},\emptyset)$,
where $J= (0,1)$ or $[0,1)$ or $[0,1],$ $\hat F$ has an ideal triangulation.
\ecor



\begin{figure}[h]
\centerline{\diag{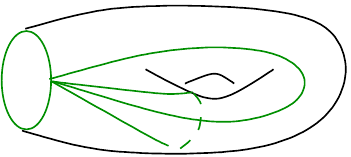}{1.8in}} 
\caption{An ideal triangulation $\Gamma=\{\gamma_0,\gamma_1,\gamma_2\}$ (in green) of a torus with one boundary component.}
\end{figure}

The proof of finite generation of the reduced Kauffman bracket skein algebras of punctured surfaces in \cite{AF, PS2} uses the fact that the reduced Kauffman bracket skein algebra has a basis of multi-curves, each of which has a canonical position with respect to any ideal triangulation of $\hat F.$ Specifically, each multi-curve in $(F,M)$ can be placed so that it consists of mutually disjoint simple (unoriented) arcs in each face of the ideal triangulation.
Such canonical position is unique. As alluded earlier, we prove an analogous statement for web diagrams in Sec. \ref{s-joy}.

%
\section{The Graded Skein Algebra}
\lb{s_graded}
%

In this section we will consider gradings on $\cal{RS}(\hat F)$ induced by ideal triangulations $\Gamma=\{\gamma_1,...,\gamma_N\}$ of $\hat F.$ 
A \underline{weight} $w(W)$ of a web $W$ in $\hat F$ with respect to $\Gamma$ (as above) is the minimal geometric intersection number of $W$ with $\Gamma$. We do not allow for the interior vertices and crossings of $W$ to lie in the edges $\gamma_1,...,\gamma_N$ (although we do allow isotopies moving vertices and crossings through $\gamma$'s).

Let $F_\ell$ be the subspace of $\cal{RS}(\hat F)$ spanned by webs of weight $\leq \ell.$
(Similar filtrations were considered in \cite{CM, Mu, Le1, AF, PS2}.)
Note that $\{F_\ell\}_{\ell\geq 0}$ is an algebra filtration of $\cal{RS}(\hat F)$ and, therefore,
it defines the \underline{graded reduced skein algebra}, $\cal{GRS}(\hat F)$. (Clearly, it depends on $\Gamma$ but we suppress that in the notation.) It is easy to see that $\cal{GRS}(\hat F)$ has the same basis of Proposition \ref{pro-basis-RS} (or Corollary \ref{cor-basis-RS}) as $\cal{RS}(\hat F)$ has.

\bthm[Proof in Sec. \ref{s_main_proof}]
\lb{main-gr}
${\cal GRS}(\hat F)$ is finitely generated as an $R$-algebra for every finite ideal triangulation of a marked surface $\hat F.$
\ethm

\noindent{\bf Proof of Theorem \ref{t-main}:}
For a closed surface $F$ there is a natural epimorphism 
${\cal RS}(F',\emptyset)\to {\cal RS}(F,\emptyset)$
where $F'$ is $F$ with some points removed. Therefore, for the sake of proof of Theorem \ref{t-main}, it is enough to assume that $F$ is not closed.

By a further removal of points from $\hat F$ if necessary, we can ensure that $\hat F$ has a finite ideal triangulation (cf. Corollary \ref{c-triangulable}). Since any element of the associated graded algebra $GA$ of an algebra $A$ lifts to an element of $A$ and any generating set of $GA$ lifts to a generating set of $A$, Theorem \ref{main-gr} immediately implies Theorem \ref{t-main}.
\qed

The proof of Theorem \ref{main-gr} requires several ingredients discussed in the remainder of the paper. 

%
\section{Minimizing geometric intersection numbers}
\lb{s-min}
%

From now now, when speaking of a web $W$ in $\hat F$ in presence of infinite arcs, we always assume that $W$ is in general position with respect to them.

\bpro[Proof in Sec. \ref{s-indices}]\lb{p-int-red}
For any infinite arc $\gamma$ in $\hat F$, 
and any non-elliptic web $W$ in $\hat F$, any planarly isotopic image $W'$ of $W$ with a minimal geometric intersection number with $\gamma$ can be obtained from $W$ by cap, vertex, crossbar 
moves (pictured below), the flip move (cf. Fig. \ref{f-isotopic}) and an isotopy away from $\gamma$.
Note that cap, vertex, crossbar moves reduce the intersection number of $W$ with $\gamma$ and the crossbar 
move does not change it. Consequently, the above moves monotonically reduce the intersection number of $W$ with $\gamma$.
\epro

We call  cap and vertex moves \underline{intersection reduction moves}. 

\begin{figure}[h]
\centerline{\diag{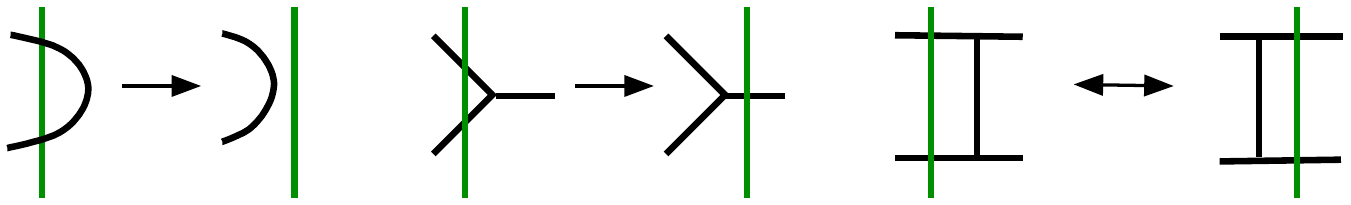}{3.5in}} 
\caption{Intersection reduction moves (cap move, vertex move) and a crossbar move. The arc $\gamma$ of the triangulation is vertical green.}
\lb{f-reduction-moves}
\end{figure}

\bcor\lb{c-min}
Every crossingless web in $\hat F$ can be isotoped though cap, vertex, and crossbar moves so that it minimizes its geometric intersection number with each arc of any disjoint collection of infinite 
arcs of $\hat F$ simultaneously.
\ecor

\bpr
It is an immediate consequence of the fact that (a) none of the moves of Proposition \ref{p-int-red} increases the intersection numbers with other infinite arcs and (b) the flip move is independent of others (i.e. it commutes with others) and does not change the intersection numbers with other infinite arcs.
Therefore, it is unnecessary for achieving a minimal intersection number.
\epr

A web as in Corollary \ref{c-min} is said to be in a \underline{minimal position} with respect to 
a collection of infinite arcs. 

\bpro\lb{p-min-pos}
A minimal position of any non-elliptic web with respect to any disjoint in $\hat F$ collection of infinite arcs is unique up to crossbar moves and flip moves (and planar isotopy within the complement of these arcs).
\epro

\bpr  The proof is by induction on the number $n$ of infinite arcs $\gamma_1,...,\gamma_n$. 
For $n=1$ it follows from Proposition \ref{p-int-red} and the fact that only the crossbar and flip moves can be applied to a diagram in a minimal position. For the inductive step, suppose that the statement holds for $n$ and suppose $W$ and $W'$ are two isotopic reduced webs, both in minimal position with respect to infinite arcs $\gamma_1,...,\gamma_{n+1}$. Since $W,W'$ are minimal with respect to $\gamma_{n+1},$ They can be made isotopic in $\hat F-\gamma_{n+1}$ by flips and crossbar moves.
Now the statement follows from the inductive assumption.
\epr

%
\section{Face indices of surface graphs. Proof of Prop. \ref{p-int-red}}
\lb{s-indices}
%

A \underline{surface graph} $G$ in $F$ is a compact graph properly embedded in $F$. Zero-valent vertices are not allowed and the monovalent vertices of $G$ must coincide with $G\cap \p F$. We allow for loop components in $G$ though.

Closures (in $F$) of connected components of $F-G$ are called \underline{faces} of $G$. 

Our proof of Proposition \ref{p-int-red} relies on the notion of index of a face. Specifically, the index of a face $f$ of $G$ is
\beq 
ind(f)=\chi(f)-\frac{int(f)}{2}-ext(f)+\sum_{v} \frac{1}{val(v)},
\eeq
where
\bite
\item $\chi(f)$ is the Euler characteristic of $f$ 
\item $int(f)$ and $ext(f)$ are the numbers of interior and exterior edges of $f$. (An exterior edge is a segment of $\p F$ which may or may not be an edge of $G$.) 
\item $val(v)$ is the number of faces a vertex $v$ of $G$ belongs to.
\eite 
If an interval appears twice as an edge of a face $f$ then it is counted twice in $int(f).$ Similarly, if a vertex appears many times in $\p f$, $val(v)$ counts that face with a multiplicity.
(Note that components of $\p f$ which are internal loops of $G$ or coincide with components of $\p F$, do not contribute to that sum. Ends of $F$ do not contribute to it either.)


\blem
\beq\lb{e-E-i}
\chi(F)=\sum_f ind(f) 
\eeq
where the sum is over all faces $f$ of $G.$
\elem

We will need it only for surface graphs $G$ in compact $F$ with contractible faces. In this case this equality follows from the fact that $G$ defines a cell decomposition of $F$, and since the Euler characteristics of every face is $1$, $\sum_f ind(f)$ adds up to the number of $2$-cells minus the number of edges ($1$-cells) plus the number of vertices ($0$-cells).

From now on we will assume that $G$ is a web in $F$ (with some possible endpoints in $\p F$). In particular, all internal vertices of $G$ are $3$-valent.

\begin{figure}[h] 
\begin{center} \diag{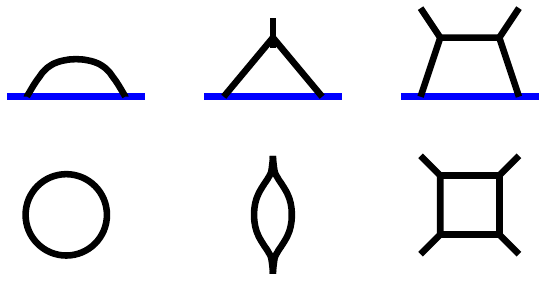}{1.5in}
\end{center}
\caption{Positive index faces. Blue line denotes $\p F$.}
\label{f-pos} 
\end{figure}

Figure \ref{f-pos} shows all possible positive index faces for surface graphs. 
From left to right, the \underline{exterior bigon} has index $\frac{1}{2}$, the \underline{exterior trigon} has index $\frac{1}{3}$, the \underline{exterior quadrigon} has index $\frac{1}{6}$, the disk has index $1$, the \underline{monogon}  has index $\frac{1}{2}$, the \underline{bigon} has index $\frac{2}{3}$, and the \underline{quadrigon} has index $\frac{1}{3}$.

The proof of Proposition \ref{p-int-red} is based on:

\blem\lb{l-bigon-red} 
Any non-elliptic web $W$ in $\hat F$ in general position with respect to disjoint, parallel infinite arcs $\gamma_1,\gamma_2$ bounding a bigon $B$ can be isotoped through the intersection reduction moves (through $\gamma_1$ and $\gamma_2$) and crossbar moves so that it intersects $B$ in a collection of disjoint, parallel arcs connecting $\gamma_1$ with $\gamma_2$, depicted by vertical arcs in figure below. 
\elem
\begin{center}
\begin{picture}(2,.3)
\put(0,0){\diag{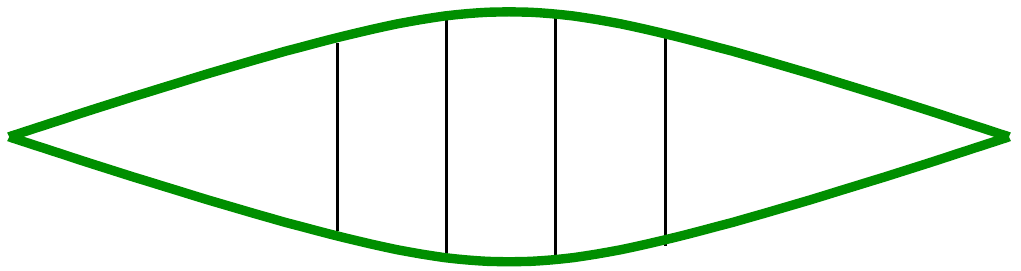}{1.8in}}
\put(1.5,-0.22){$\gamma_1$}
\put(1.5,0.2){$\gamma_2$}
\put(-0.2,0){$m_1$}
\put(1.85,0){$m_2$}
\end{picture}
\end{center}\vspace*{.2in}

\noindent{\it Proof by contradiction:} Suppose that the statement is false. Among webs $W$ contradicting the statement consider a web $W$ with the minimal number of vertices in $B$. Let $W'$ be $W\cap B$ with all 
simple arcs connecting $\gamma_1$ and $\gamma_2$ (pictured as vertical segments in the picture above) removed.
Formula (\ref{e-E-i}) implies that the sum of indices of faces of $W'$ in $B$ is $1$, the Euler characteristic of $B$. Let $f_1,f_2$ be the faces of $W'$ containing $m_1$ and $m_2$.  By discussion above (cf. Figure \ref{f-pos}) $f_1,f_2$ have indices $\leq 1/2.$ Furthermore, none of them can be $1/2$ since that would imply that $f_1$ or $f_2$ is a bigon and, hence, $W'$ contains a simple arc -- contradicting the definition of $W'.$

Hence, $W'$ must contain another face of positive index in $B$. Since $W'$ is non-elliptic, that face must be an external bigon, trigon, or quadrigon. That face can be eliminated by either an intersection reduction move (cap or vertex move) or a crossbar move (away from $m_1,m_2$) contradicting the minimality of vertices of $W\cap B.$
\qed\vspace*{.1in}

\noindent{\bf Proof of Proposition \ref{p-int-red}:} By applying a sequence of cap, vertex, and crossbar moves to $W$ we can assume that no further reduction of $|W\cap \gamma|$ through those moves is possible. The proof continues by contradiction: Assume that an arc $\gamma'$ is isotopic to $\gamma$ rel its endpoints with 
\beq\lb{p-contr-ass}
|\gamma' \cap W| < |\gamma \cap W|.
\eeq 
By applying an arbitrarily small isotopy to $\gamma$ and $\gamma'$ if necessary, one can assume that $\gamma$ and $\gamma'$ are in general position and that $|\gamma\cap \gamma'|$ is finite. 
Let us assume for $\gamma, \gamma'$ as above that $|\gamma\cap \gamma'|$ is as small as possible. By a theorem of Epstein, \cite{Ep}, there is a bigon $B$ bounded by a subarc $\alpha$ of $\gamma$ and a subarc $\alpha'$ of $\gamma'.$ 

Finally, let as assume that for $\gamma,\gamma',W$ and $B$ as above, the number of vertices of $W$ in $B$ is as small as possible. Given these assumptions, Lemma \ref{l-bigon-red} (applied to
$\hat F$ with the endpoints of $\alpha,\alpha'$ removed) implies that $W$ intersects $B$ in a collection of simple, parallel, arcs connecting $\alpha$ with $\alpha'$.

If $\gamma\cap \gamma'= \emptyset$ then $\alpha=\beta,$ $\alpha'=\beta'$ and the above statement contradicts
(\ref{p-contr-ass}). If $\gamma\cap \gamma'\ne \emptyset$ then pushing $\gamma$ through $B$ reduces the number of intersections with $\gamma'$ -- contradicting the assumption about the minimality of $|\gamma\cap \gamma'|$. 
\qed

%
\section{Webs in bigons}
\lb{s-bigons}
%
 
Let $B$ be the (ideal) bigon $D^2-\{m_1,m_2\}$  with its boundary intervals (the connected components of $\p D^2-\{m_1,m_2\}$) denoted by $\p_+ B,$ $\p_- B$, oriented from $m_1$ to $m_2,$ and so that $\p_+ B$ lies on the left and $\p_- B$ on the right assuming (as usual) the counterclockwise orientation of $B$. 

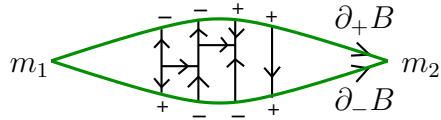
\begin{figure}[h]
\begin{center}
\begin{picture}(2,0.5)
\put(0,0.1){\diag{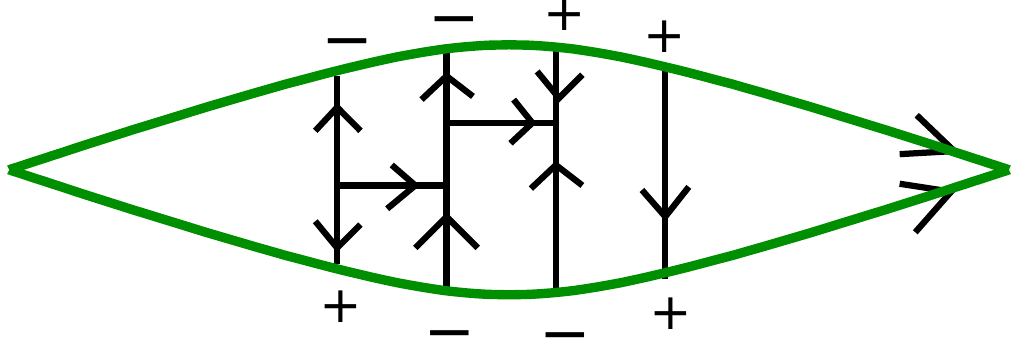}{1.8in}} 
\put(1.5,-0.1){$\p_- B$}
\put(1.5,0.33){$\p_+ B$}
\put(-0.2,0.1){$m_1$}
\put(1.85,0.1){$m_2$}
\end{picture}
\end{center}
\caption{A crossbar web in bigon $B$ with plus signature $(-,-,+,+)$ and minus signature $(+,-,-,+)$.}
\lb{f-bigon-conventions}
\end{figure}


A \underline{crossbar web} of \underline{index} $n$ in $B$ consists of $n$ \underline{vertical lines}, which are simple, mutually disjoint arcs connecting $\p_+ B$ and $\p_- B$ and a certain number of horizontal intervals (``crossbars") connecting some of the adjacent vertical lines, as in Figure \ref{f-bigon-conventions}. (Together they form a graph, which one orients according to web orientation rules.) 
In particular, there are $2^n$ index $n$ crossbar webs without crossbars, called \underline{trivial}.
(Our crossbar webs are the ladder webs of \cite{CKM}, with the only small, but inessential difference, that our crossbars are horizontal, unlike theirs.)




The signs of intersections of the external edges of a crossbar web $W$ with $\p_+ B$ and with $\p_-B$  
form the \underline{plus-signature} $(p_1,...,p_n)\in \{\pm\}^s$ and the  \underline{minus-signature} $(q_1,...,q_n)\in \{\pm\}^s$, of $W$, respectively, cf. Figure \ref{f-bigon-conventions}.



We are going to consider \underline{braid diagrams} drawn in the bigon in $B$, connecting $\p_+ B$ with $\p_- B.$ 
However, unlike for typical braid diagrams, their strands may have arbitrary orientation.
We say that a braid diagram is \underline{proper} if no two strings of it with coinciding orientations cross. 
(Note that we do not consider braid diagrams up to braid relations and that our notion of proper braid diagram is non-standard.) The signs partition $\{1,...,n\}$ into two sets and a proper braid is a shuffle with respect of these sets -- that is it a permutation of $\{1,...,n\}$ preserving the interior order of these two sets.


A proper braid diagram in $B$ is \underline{minimal} if no two strands in it cross more than once.

\blem\lb{l-min-cons}
A proper braid diagram in $B$ is minimal iff no two strands in it cross twice consecutively (forming a bigon).
\elem

\bpr The implication $\Rightarrow$ is obvious. Proof of $\Leftarrow$: Let $\beta$ be a braid diagram in which two strands cross (at least) twice. The arcs connecting these two crossings form a bigon, which $\beta$ may intersect. There may be a number of bigons of this form cut out by $\beta$. Let $\cal B$ be a minimal one -- that is one which does not contain a smaller one inside.
Since $\beta$ is proper, no strand of $\beta$ may intersect both sides of $\cal B$. Therefore, either
\bite
\item $\beta$ does not intersect the interior of $\cal B$, implying that $\cal B$ is formed by two strands in it crossing twice consecutively, or
\item $\beta$ intersects $\cal B$ in an arc with both endpoints on one side of $\cal B$. That double intersection defines a smaller bigon in $\cal B$ -- contradicting the minimality of $\cal B.$
\eite
\epr




Note that the I-resolution of crossings (cf. Fig \ref{fig-cross-res}) in a proper braid diagram $\beta$ yields a crossbar web. 
We denote it by $H_\beta.$ Furthermore, every crossbar web is obtained that way.
(Note however that (a) isotopic braid diagrams may define different crossbar webs and (b) $H_\beta$ is insensitive to crossing changes in $\beta$.) 
$H_\beta$, for $\beta$ minimal, is called a \underline{minimal crossbar web}.

For the sake arguments in the remainder of this paper it is useful to consider the notion of \underline{weakly reduced webs} by which you mean non-elliptic webs without exterior bigons and triangles (the first two diagrams in Fig. \ref{fig-degeneracies}), but which may have triple edge paths parallel to $\p F.$ 
Clearly, every minimal crossbar web in $B$ is weakly reduced. Furthermore, we have:

\blem\lb{l-red-signature}
(1) Every weakly reduced web $W$ is a minimal crossbar web.\\ 
(2) A minimal $\beta$ such that $W=H_\beta$ is unique up to crossing sign changes.\\
(3) Every weakly reduced web is determined by its signature uniquely (in the set of crossbar webs).
\elem

\bpr
(1) By Lemma \ref{l-bigon-red}, every web in a bigon can be trivialized by intersection reduction and crossbar moves. Since a weakly reduced web does not allow for intersection reduction, it must be a crossbar web,
$W=H_\beta$ for some proper $\beta$. If $\beta$ is not minimal then it contains two strings with consecutive crossings by Lemma \ref{l-min-cons}. These consecutive crossings 
yield crossbars which together with vertical lines connecting them form a $4$-gon in $W$, contradicting the non-ellipticity of $W$.

(2) and (3) follow from the fact that (a) the signature of $\beta$ coincides with that of $W=H_\beta$ and (b) this signature determines a minimal braid diagram uniquely.
\epr

%
\section{Canonical webs in triangles}\lb{s-triangle}
%

A disk with three marked points in its boundary is called a \underline{triangle}. 
An arbitrarily chosen side of a triangle $T$ will be drawn horizontally and called the \underline{base}. 
The \underline{degree} $d\in \Z_{>0}$ \underline{pyramid} $H_d$  in $T$ is a web consisting of: 
\bite
\item  $d$ horizontal lines parallel to the base, which we number from $1$ (the shortest) to $d$ (the longest)
\item $i$ vertical intervals connecting the $i$-th line with the $i+1$-th line, for $i=1,2,...,d$, so that all vertical lines are in a staggered pattern. Here, the $(d+1)$-th line is the base (which is not the part of $H_d$).
These vertical intervals split the horizontal lines into horizontal intervals.
\item the horizontal and vertical intervals are oriented according to the web orientation rules, so that all endpoints of $H_d$ point outwards,
\eite
c.f. Figure \ref{f-pyramids}. 

\begin{figure}[h]
\centerline{\diag{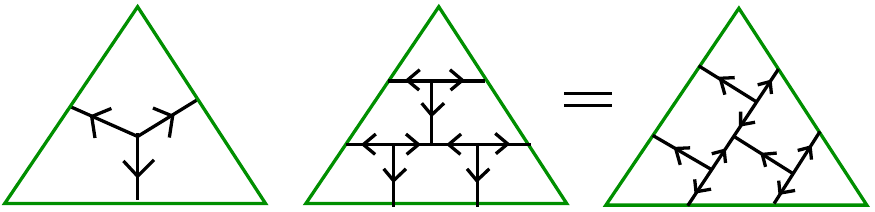}{2.5in}}
\centerline{\diag{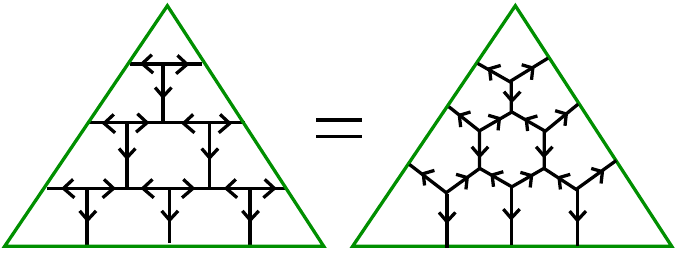}{2.3in}}
\caption{Top: A pyramid of degree $1$ and a pyramid of degree $2$ with two different choices of the base. Bottom: A pyramid of degree $3$. Note the $2\pi/3$ rotational symmetry of the honeycomb form of the pyramid.}
\lb{f-pyramids}
\end{figure}

Note that $H_d$ does not depend on the choice of base of $T$, cf. Figure \ref{f-pyramids}.
By reversing all orientations of $H_d$ we obtain $H_{-d}.$ Therefore, $H_d$ is defined for all $d\in \Z$ with $H_0=\emptyset.$

We say that a web in $T$ is \underline{canonical} if it consists of $H_d$, for some $d\in \Z$,
and of a number of simple oriented arcs, mutually disjoint from each other and from $H_d,$ each connecting different sides of the ideal triangle,  cf. Figure \ref{f-canonical}. These are \underline{corner arcs}.

\begin{figure}[h]
\centerline{\diag{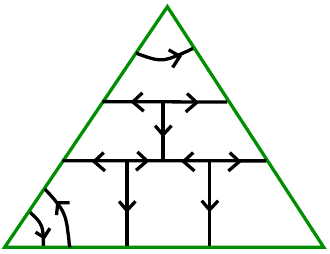}{1.5in}}
\caption{An example of a canonical web in a triangle.}
\lb{f-canonical}
\end{figure}

%
\section{Joy-sparking and canonical webs}\lb{s-joy}
%

Let $\hat F$ be a marked surface with an ideal triangulation $\Gamma=\{\gamma_1,...,\gamma_N\}$.
Choose arbitrary orientations of $\gamma_1,...,\gamma_N$ and their disjoint tubular neighborhoods $\cal N(\gamma_1), ..., \cal N(\gamma_N).$ 

Recall from Sec. \ref{s-ideal}, that all boundary arcs are included among the edges of $\Gamma$; we call them \underline{external}. 
Each external edge $\gamma_i$ necessarily coincides with one of the two boundary intervals of $\cal N(\gamma_i).$ 


We call $\cal N(\gamma_1), ..., \cal N(\gamma_N)$ a \underline{padded ideal triangulation} of $\hat F$ and $\cal N(\gamma_1), ..., \cal N(\gamma_N)$ are called its \underline{bigons}. Each of them is an ideal bigon. (In the case of $\hat F=S^1\times (0,1)$, it is a bigon with both of its sides identified.)

The connected components of $F-\bigcup_{i=1}^N \cal N(\gamma_i)$ are its \underline{faces}.
By the definition of an ideal triangulation, each of them is an ideal triangle.
 
The orientations of $\gamma_1,...,\gamma_N$ define the oriented boundary intervals $\p_+ \cal N(\gamma_i),$ $\p_- \cal N(\gamma_i)$ of the bigons $\cal N(\gamma_i)$ for $i=1,...,N,$ as in Sec. \ref{s-bigons}. 

From now on we consider non-elliptic webs only, unless specifically stated otherwise.
Taking advantage of Corollary \ref{c-min}, we will assume that these webs are in the minimal position with respect to infinite arcs $\p_\pm \cal N(\gamma_1), ..., \p_\pm \cal N(\gamma_N)$ and simply say that they are in a \underline{minimal position}.

By Proposition \ref{p-min-pos}, such position is unique up to crossbar and flip moves.
We call crossbar moves pushing crossbars into bigons \underline{tidying up} operations, 
cf. Figure \ref{f-tidy-up}.
\begin{figure}[h]
\centerline{\diag{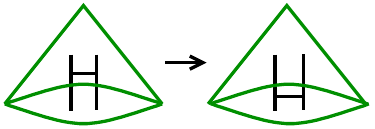}{2in}}
\caption{Tidying up operation.}
\lb{f-tidy-up}
\end{figure}

Consider a non-elliptic web $W$ in minimal position and apply to it tidying up operations.
Inspired by Marie Kondo, \cite{Ko}, we say that such web \underline{sparks joy} if no further tidying up operations are possible.
As a consequence of Proposition \ref{p-min-pos} we have:

\bcor\lb{c-unique-joy}
The result of putting a non-elliptic web $W$ (considered up to planar isotopy) in a joy-sparking position (with respect to 
$\cal N(\gamma_1),...,\cal N(\gamma_N)$) is unique up to flip moves and crossbar passes (cf. Fig. \ref{f-crossbar-pass}). 
\ecor

\begin{figure}[h]
\centerline{\diag{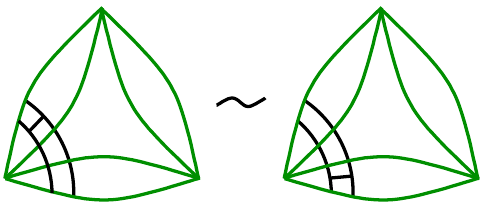}{1.8in}}
\caption{A crossbar pass. Orientations not specified.}
\lb{f-crossbar-pass}
\end{figure}

In order to characterize joy-sparking webs, it is best to narrow focus to the weakly reduced webs introduced in Sec. 
\ref{s-bigons}. We say that a weakly reduced web $W$ in $\hat F$ is in a \underline{canonical position} with respect to an ideal triangulation $\Gamma$ if $W\cap T$ is canonical for every face $T$ of the triangulation and
$W\cap \cal N(\gamma_i)$ is a minimal crossbar web for every $i$.

\bthm\lb{t-main-canonical}
A weakly reduced web is in a minimal position and sparks joy iff it is in a canonical position.
\ethm

\noindent{\bf Proof of $\Leftarrow$:}
A web in a canonical position is in a minimal position and sparks joy, since it does not allow for any tidying up operations. 

\noindent{\it Proof of $\Rightarrow$:} Since $W$ is in a minimal position with respect to $\p_{\pm} \cal N(\gamma_i)$, Lemma \ref{l-bigon-red} implies that $W\cap \cal N(\gamma_i)$ is a crossbar web. It is minimal because $W$ is weakly reduced. Therefore, it is enough to prove that  $W$ is in canonical position in every triangle $T$ of the triangulation. By Lemma \ref{l-bigon-red}, the only connected components of $W$ which do not touch all sides of $T$ are corner arcs. Let $W'$ be $W$ stripped of all these corner arcs. 

Let us assume that $W'\ne \emptyset$ now since otherwise we are done. Consider the faces $f_0,f_1,f_2$ of $W'$ in $T$, as in Sec. \ref{s-indices}, containing the vertices $v_0,v_1$ and $v_2$ of $T$, respectively. By the above discussion, $f_i$ is bounded by a path of $n_i$ edges of $W'$ for $n_i\geq 2.$ Since $f_i$ has one external edge (with $v_i$ inside of it),
$$ind(f_i)=1- n_i/2-1+2\cdot \frac{1}{2}+(n_i-1)\frac{1}{3}=\frac{1}{3}-\frac{n_i-2}{6}.$$
Since $W'$ is non-elliptic web without $1$-, $2$- and $3$-gons parallel to the three sides of $T$, the faces $f_0,f_1, f_2$ are the only faces with potentially positive index. Since the sum of the face indices is $\chi(T)=1,$ 
we see that $ind(f_i)=\frac{1}{3}$ and $n_i=2$ for $i=0,1,2.$ That implies that all faces other than $f_0,f_1,f_2$ have index $0$. In particular, all internal faces of $W'$ are hexagonal and all external faces other than $f_0,f_1,f_2$ are bounded by a 4-gon, as in Fig. \ref{f-pyramids}. Consequently, $W'$ is a pyramid. 
\qed\vspace*{.1in}

%
\section{Web coordinates}\lb{s-coordinates}
%

The \underline{intersection coordinates} 
of a web $W$ in a canonical position with respect to a triangulation $\Gamma$ are given by the numbers of $+$ and $-$ signs of the intersections of $W$ with the bigons of the triangulation.
Specifically, the intersection coordinates of $W$ are 
\bite
\item the numbers, $e_{+,i}(W),$ of positive endpoints of $W$ in $\p_+ \cal N(\gamma_i)$ or, equivalently, in $\p_- \cal N(\gamma_i)$. (Because $W\cap B$ is a crossbar web, these numbers of positive endpoints coincide, 
cf. Sec. \ref{s-bigons}.)
\item the numbers, $e_{-,i}(W),$ of negative endpoints of $W$ in $\p_+ \cal N(\gamma_i)$ or, equivalently, in $\p_- \cal N(\gamma_i)$.
\eite
The coordinates $e_{\pm,i}(W)$ are indexed by $i=1,...,N,$ where $N$ is the number of edges in $\Gamma.$

Note that these coordinates alone do not determine the isotopy class of a canonical web in a triangle, cf. 
Fig. \ref{fig-different-winding}.
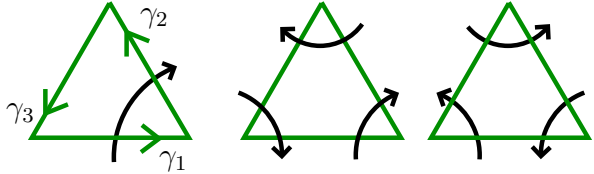
\begin{figure}[h]
\begin{center}
\begin{picture}(3,0.8)
\put(0,0.35){\diag{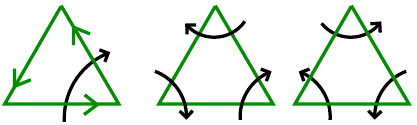}{3in}} 
\put(.7,-.05){$\gamma_1$}
\put(.6,0.7){$\gamma_2$}
\put(-0.1,0.2){$\gamma_3$}
\end{picture}
\end{center}
\caption{Left: A fragment of a web with coordinates $e_{-,1}=e_{+,2}=r=1$ and $e_{+,1}=e_{-,2}=e_{\pm,3}=0$. Right: Non-isotopic canonical webs with coinciding intersection coordinates. Their rotation numbers are $3$ and $-3$.}
\lb{fig-different-winding}
\end{figure}
For that reason we need the \underline{rotation number} of a web $W$ in canonical position in a face $T$,
defined as
$$r_T=\sum_\alpha \ve(\alpha),$$ 
where the sum is over all corner arcs of $W\cap T$ and
$\ve(\alpha)=1$ if $\alpha$ is oriented clockwise and $\ve(\alpha)=-1$ otherwise.

The two intersection coordinates per each edge $\gamma_i$, and the rotation number per each face of the triangulation are simply called the \underline{coordinates} of a web in a canonical position. 

According to Corollary \ref{c-unique-joy}, the joy-sparking position of a weakly reduced web is unique up to flip moves and up to crossbar passes. Since these operations do not affect the coordinates, we say that the coordinates of any weakly reduced web $W$ are those of $W$ put into a canonical position.
By Sec. \ref{s-reduced-skein} that two weakly reduced webs are equivalent if they are related by isotopy and by permutations of parallel arc components. 
Since these operations preserve coordinates, we obtain:

\bcor\lb{c-coord-well-def}
The above coordinates are well defined on weakly reduced webs considered up to equivalence. 
\ecor

These coordinates do not determine weakly reduced webs in when $\p \hat F\ne \emptyset,$ because a boundary bigon may contain an arbitrary minimal crossbar web and those are not determined by intersection coordinates.
However, somewhat surprisingly we have:

\bthm[Proof in Sec. \ref{s-coord-proof}]\lb{t-coordinates}\ \\
Each reduced web in $\hat F$ is determined up to equivalence by its coordinates. 
\ethm

As alluded in the introduction we believe that the above coordinate system leads the notion of quantum $SL(3)$-trace which embeds our reduced skein algebras into quantum tori. In particular, we propose

\bcon\label{c-domain}
For every marked surface $\hat F$ and any domain $R$,
${\cal RS}(\hat F)$ is an Ore domain. 
\econ

(This statement seems to be a consequence of the recent \cite{Ki}.)

For a web $W$ in a canonical position, the degree of the pyramid in $T\cap W$ is simply called the
\underline{degree} of $W$ in $T$. The proof of Theorem \ref{t-coordinates} will use the following observation: If $T$ is bounded by $\gamma_{i_0}, \gamma_{i_1},$ and $\gamma_{i_2},$ with the orientations induced by that of $T$, then the degree of $W$ in $T$ is
\beq\lb{e-honeycomb}
\frac{1}{3}\sum_{k=0}^2 e_{+,\gamma_{i_k}}(W)-e_{-,\gamma_{i_k}}(W).
\eeq

If $\gamma_{i_k}$ is oriented in the opposite way to $\p T$ then $e_{+,\gamma_{i_k}}(W)$ and 
$e_{-,\gamma_{i_k}}(W)$ need to be interchanged in the above formula.


%
\section{Relation to 
Fock-Goncharov-Shen theory}
\lb{s-FG}
%

Recall that a flag in $\C^n$ is a maximal sequence of vector subspaces 
$$V_0\subsetneq V_1\subsetneq ... \subsetneq V_{n-1}\subsetneq V_n=\C^n$$ 
(hence, $dim_\C V_k=k$) and an \underline{(affine) $SL(n,\C)$-flag} is a flag together with a choice of a vector
$v_i\in V_i/V_{i-1}$ for every $i=1,...,n,$ such that $det(v_1,...,v_n)=1.$ 
Obviously, every affine $SL(n,\C)$-flag is determined by the $n$-tuple $v_1,...,v_n$ as above.
The group $SL(n,\C)$ acts transitively on the space of affine flags on the right and the stabilizer of the standard flag, with $V_k=\ \text{Span}\{e_1,...,e_k\},$ for $k=1,...,n,$ is the group of upper triangular matrices with ones on the diagonal, denoted by $U$. (Note that $U$ is a maximal unipotent subgroup of $SL(n,\C).$) 
Consequently, the space of affine $SL(n,\C)$-flags can be identified with $SL(n,\C)/U.$

Let $G$ be a reductive group with a maximal unipotent group $U\subset G.$
Each $G$-local system $\cal L$  on a surface $F$ (with the action of $G$ on the right) defines the associated
\underline{principal affine bundle} on $\hat F$,
\beq\label{e-principal-U}
 \cal L_{\cal A} =\cal L\times_{G} U.
 \eeq

A \underline{decorated $G$-local system} on $\hat F$ is a pair $(\cal L, \alpha)$, where $\cal L$ is a $G$-local system on $F$, and $\alpha$ is a flat section of the restriction of $\cal L_{\cal A}$ to $\p \hat F$. 
An \underline{essentially decorated $G$-local system} on $\hat F$ is a pair $(\cal L, \alpha)$, where $\cal L$ is a $G$-local system on $F$, and $\alpha$ a flat section of the restriction of $\cal L_{\cal A}$ to the boundary intervals of $\p \hat F$. (Hence, essentially decorated local systems do not contain decoration for unmarked boundary components of $S$.)
Following \cite{FoGo, GS1}, we denote the spaces of decorated and essentially decorated twisted $G$-local systems on $\hat F$ by $\cal A_{G, \hat F}$ and by $\text{Loc}_{G,\hat F}$ respectively. (We do not define twisting here, since it can be ignored for $SL(3)$.)
Both are affine varieties. 

As observed in \cite{FP}, there is a bijection
\beq\label{e-correspondence}
\text{
$SL(3,\C)$-flags\,$(v_1,v_2,v_3)$ $\leftrightarrow$ pairs\,$(v,f)\!\in\!(V\!-\!\{0\})\otimes (V^*\!-\!\{0\})$}
\eeq 
such that $f(v)=0.$ In this correspondence, $v=v_1$ and $f$ is such that $f(v_1)=f(v_2)=0$ and 
$f(v_3)=1.$ 

With the aid of the above correspondence, every $3$-web $W$ in $\hat F$ defines a function on $\text{Loc}_{SL(3), \hat F}(\C)$ as follows: 
for each $\cal L$ on $\hat F$, each edge of $W$, going from $p\in F$ to $q\in \hat F$ defines linear transformation $V_p\to V_q$, where $V_p\simeq \C^3$ is the fiber over $p$ of the bundle $\cal L \times_{SL(3,\C)}\C^3$. A sink at $p$ represents the volume form
$V_p^{\otimes 3}\to \wedge^3 V_p=\C$ and a source at $p$ represents the dual volume form.
To an endpoint of $W$ in $\p \hat F$ one associates either the vector $v$ or the dual vector $f$ of the pair $(v,f)$
representing the section of $\cal L_{\cal A}$ at that point. 
A contraction of all these tensors for a given web $\Gamma$ defines a function on $\text{Loc}_{SL(3), \hat F}(\C)$.
This is the construction of  \cite{FP}, except that theirs is related to ours by the composition with the natural projection 
$${\cal A}_{SL(3), \hat F}(\C)\to \text{Loc}_{SL(3), \hat F}(\C).$$
(Similar versions of this construction appear for example in \cite{S1} and \cite{Ku}.)

Webs considered as functions as above, satisfy the skein relations of the reduced skein algebra, \eqref{skein1}-\eqref{skein3} and \eqref{e-bskein1}-\eqref{e-bskein2} for $q=a=1.$
Therefore, the above construction factors to a homomorphism 
$$\phi_3: {\cal RS}(\hat F)_{q=a=1}\to \cal O(\text{Loc}_{SL(3), \hat F}(\C)),$$
where the ring of coefficients is $R=\C.$ We call it the \underline{tensor contraction} \underline{homomorphism}.
Note that our $q=a=1$ reduced skein algebra relations coincide with those of \cite{FP}, except for the additional relation \diag{brel4}{.6in} $= 0$ in ${\cal RS}(\hat F).$

Now \cite[Conj. 5]{FP} can be reformulated and strengthened to:
\bcon\label{c-phi3}
$\phi_3$ is an isomorphism.
\econ

Fock-Goncharov proved that $\cal A_{G,\hat F}$ has an atlas with transition functions subtraction free. Such an atlas defines a tropical version of this variety, $\cal A_{G,\hat F}(\Z^t)$. 
Goncharov-Shen, \cite{GS1}, defined a rational function on $\cal A_{PGL(n),\hat F}(\C)$, called a potential, which tropicalizes to 
$$\cal W: \cal A_{PGL(n),\hat F}(\Z^t)\to \Z.$$ 
The preimage $\cal W^{-1}(\Z_{\geq 0})$ is called the set of positive points and denoted by $\cal A_{PGL(n),\hat F}^+(\Z^t)$.

Goncharov-Shen conjectured a version of a mirror duality, relating $\cal A_{PGL(n),F}^+(\Z^t)$ and $\cal O(\text{Loc}_{SL(n), \hat F}).$ They proved in fact (\cite[Thm. 10.12]{GS1}) that every triangulation of $\hat F$ defines a basis of $\cal O(\text{Loc}_{SL(n), \hat F})$ parametrized by the points of $\cal A_{PGL(n),\hat F}^+(\Z^t)$. Furthermore, they conjectured the existence of a mapping-class-group-equivariant basis of $\cal O(\text{Loc}_{SL(n), \hat F})$ indexed by the points of $\cal A_{PGL(n), \hat F}^+(\Z^t)$. Our Conjecture \ref{c-phi3} implies that Corollary \ref{cor-basis-RS}
does indeed provide a mapping-class-group-equivariant basis of $\cal O(\text{Loc}_{SL(3), \hat F})$.

Their conjecture motivates:

\bcon
There is a $1$-$1$ correspondence between reduced non-elliptic webs in $\hat F$ without British highways and the points of $\cal A_{PGL(3), \hat F}^+(\Z^t)$. 
\econ

The relation between reduced non-elliptic webs and  $\cal A_{PGL(3), \hat F}^+(\Z^t)$ is investigated in further detail in \cite{DS1, DS2}. The above conjecture is also related to recent work of Kim, \cite{Ki}, in which he establishes a similar web interpretation for the points of $\cal A_{SL(3),\hat F}(\C)$.

%
%

Let us assume now that $F$ is compact and with all boundary components marked and that $\hat F$ has an ideal triangulation. Then our ${\cal RS}(\hat F)_{a=1}$ is an $SL(3)$-version of Muller's Kauffman bracket skein algebra, $Sk_q(F)$, \cite{Mu}. Quantum cluster algebra structure on ${\cal RS}(\hat F)_{a=1}$ is studied in the recent \cite{IY}.





%
\section{Proof of Theorem \ref{t-main-vr}}
\lb{s_t-main-vr-proof}
%

A \underline{size} of a web diagram is the number of its vertices plus twice the number of crossings. 
Size defines a filtration of $\cal{RS}(F,M)$ (which is not an algebra filtration).

A union of web diagrams stacked one on top of other, is called a  \underline{stack}.
Specifically, it is a web diagram of the form $W_1 \cup ... \cup W_n,$ such that $W_i$ lies on top of $W_{i+1}$ for $i=1,...,n-1$. It is like a product of webs, except that the endpoints of webs in a stack do not have to be arranged according to the rule of the multiplication of webs. We have
\beq\lb{e-boundary-cross1}
\diag{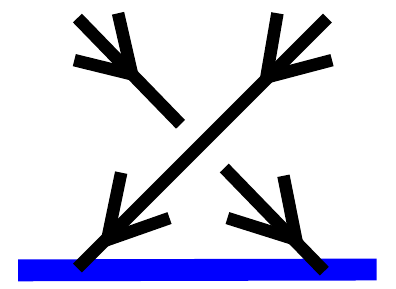}{.4in}=q^{2/3}\cdot \diag{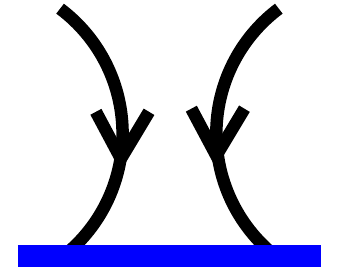}{.35in}+q^{-1/3}\cdot\diag{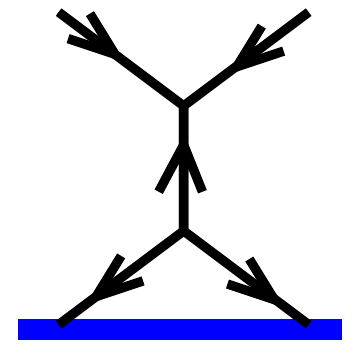}{.35in}=q^{2/3}\cdot 
\diag{smoothob}{.35in}
\eeq 
and
\beq\lb{e-boundary-cross2}
\diag{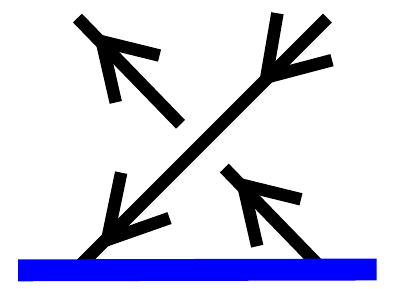}{.4in}=q^{1/3}\cdot \diag{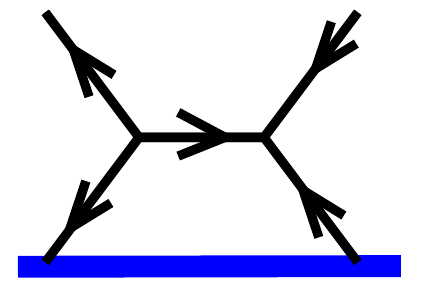}{.45in}+q^{-2/3}\cdot
\diag{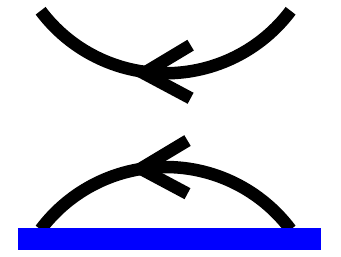}{.4in}=q^{1/3}\cdot a^{-1}\diag{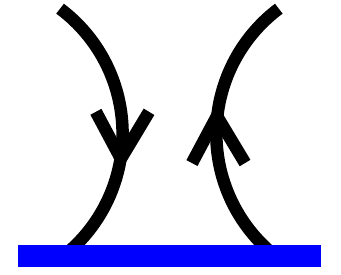}{.35in},
\eeq
where the blue horizontal line denotes $\p F$ as usual.
Hence, a stack $W_1 \cup ... \cup W_n$ equals to $W_1 \cdot ... \cdot W_n$ times a power of $q$ and $a$ in $\cal{RS}(F,M)$.

Theorem \ref{t-main-vr}(1) follows by induction on size from the lemma below. 

\blem\lb{l-stack}
Every web diagram $W$ in $(F,M)$ of size $s$ equals in ${\cal RS}(F,M)$ to a power of $q$ times a stack of descending knot, arc, and triad diagrams of total size at most $s$ plus a linear combination of webs of size less than $s$.
\elem

\bpr Let $W$ be a web in $(F,M)$ of size $s$. By (\ref{skein3}), up to a multiplicative factor of power of $q,$ one can replace all internal edges of $W$ by crossings of arbitrary sign, plus linear combinations of webs of smaller size. Since we can control the signs of these crossings, the above operations can transform $W$ into a stack of webs diagrams without internal edges  (i.e. composed of diagrams of knots, arcs, and triads) represented by a diagram of size $s$, plus a linear combination of webs of smaller size. Furthermore, these diagrams of knots, arcs, and triads can be assumed in descending position.
\epr

Proof of part (2) of Theorem \ref{t-main-vr} follows from the fact that the only descending knots, arcs and triads in $(D^2,M)$ are trivial (i.e. crossingless). For the triad that follows from the fact that any crossing between two of its edges can be undone by a twist. The possible side effect -- twisting of the framing of the edges can be eliminated by a multiplication by a power of $q$. The fact that the endpoints of arc and triad diagrams can be assumed to lie in different boundary intervals follows from \eqref{e-bskein1}.
\qed

\qed

%
\section{Proof of Theorem \ref{t-coordinates}}
\lb{s-coord-proof}
%

The coordinates of $W$ determine the degree of $W$ in every triangle of an ideal triangulation by (\ref{e-honeycomb}). Furthermore:

\blem\lb{l-corner-arc-mult}
The coordinates of a web in canonical position determine the multiplicities of all six types of corner arcs in each triangle of the triangulation.
\elem

\bpr
Let as assume as before that a face $T$ of triangulation is bounded by $\gamma_{i_0}, \gamma_{i_1},$ and $\gamma_{i_2},$ with the orientations induced by that of $T$. (The proof for other orientations is analogous.)
Let $e_{\pm,j}=e_{\pm,\gamma_{i_j}}(W)-d$ for $i=0,1,2$, for simplicity, where $d$ is the degree of the pyramid in $T$ given by (\ref{e-honeycomb}). Then $e_{\pm,0},e_{\pm,1},e_{\pm,2}$ count the ends of corner arcs in $W\cap T$. Specifically, denote the numbers of corner arcs by 
$n_{i,j},$ for $i\ne j,$ as in Figure \ref{f-corner-arcs}.

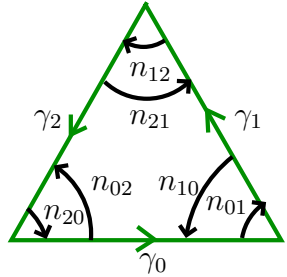
\begin{figure}[h]
\begin{center}
\begin{picture}(1.4,1.3)
\put(0,.6){\diag{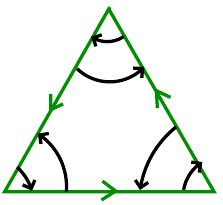}{1.5in}} 
\put(.7,-.1){$\gamma_0$}
\put(1.2,0.65){$\gamma_1$}
\put(0.15,0.65){$\gamma_2$}
\put(0.2,0.15){$n_{20}$}
\put(0.45,0.3){$n_{02}$}
\put(0.8,0.3){$n_{10}$}
\put(1.05,0.2){$n_{01}$}
\put(0.65,0.9){$n_{12}$}
\put(0.65,0.65){$n_{21}$}
\end{picture}
\end{center}
\caption{Corner arcs in a face of an ideal triangulation.}
\lb{f-corner-arcs}
\end{figure}

Then
\beq\lb{e-n2e}\left(\begin{matrix}
e_{+,0}\cr e_{-,0}\cr e_{+,1}\cr e_{-,1}\cr e_{+,2}\cr e_{-,2}\cr r 
\end{matrix}\right)=\left(\begin{matrix}
0 & 1 & 0 & 0 & 1 & 0\\
1 & 0 & 0 & 0 & 0 & 1\\
1 & 0 & 0 & 1 & 0 & 0\\
0 & 1 & 1 & 0 & 0 & 0\\
0 & 0 & 1 & 0 & 0 & 1\\
0 & 0 & 0 & 1 & 1 & 0\\
1 & -1 & 1 & -1 & 1 & -1\\
\end{matrix}\right)\cdot \left(\begin{matrix} n_{01}\\ n_{10}\\ n_{12}\\ n_{21}\\n_{20}\\ n_{02}\end{matrix}\right).
\eeq
The above $7\times 6$-matrix has zero nullity and, hence, $e_{\pm,0},e_{\pm,1},e_{\pm,2}$ and the rotation number $r$ determine the numbers $n_{ij}.$ 
\epr

\noindent{\bf Proof of Theorem \ref{t-coordinates}:} 
Let reduced webs $W_1,W_2$ in $\hat F$ have coinciding coordinates. Without loss of generality we can assume that they are in canonical position with respect to $\Gamma.$ 
By (\ref{e-honeycomb}) and Lemma \ref{l-corner-arc-mult}, $W_1$ and $W_2$ coincide in each face of the triangulation up to a permutation of parallel corner arcs in different directions.  Although the statement of the theorem is of topologically-combinatorial nature, we find it easiest to prove by applying some of the machinery of the $SU(3)$-skein algebras developed in this paper. Specifically, we are going to use the following trick: Let $a=1$ in the coefficient ring $R$. Then by Proposition \ref{pro-basis-RS} and Figure \ref{fig-parallel-arcs} two reduced webs in $\hat F$ are equivalent iff they coincide in $\cal{GRS}(\hat F).$ Therefore, we will complete the proof by showing that $W_1=W_2$ in $\cal{GRS}(\hat F).$ 

For the sake of that proof it will be convenient to expand the class of webs under consideration.
We say that a web $W$ is in an \underline{almost canonical} position with respect to a padded ideal triangulation $\cal N(\gamma_1),...., \cal N(\gamma_N)$ if it is in a canonical position in every face of the triangulation and it is a crossbar web in $\cal N(\gamma_i)$ for every $i$. However, we do not assume that $W$ is reduced and, hence, $W$ may contain $4$-gons, either internal or external.
Every internal $4$-gon involves two parallel crossbars either inside of one bigon $\cal N(\gamma_i)$ or two different ones,
cf. Figure \ref{f-4-gon-elimination}.

However, by (\ref{skein2}), the two parallel crossbars in any internal $4$-gon 
can be removed without affecting the value of the web in $\cal{GRS}(\hat F)$, cf. Figure \ref{f-4-gon-elimination}. That follows immediately from the fact that removing these crossbars corresponds to one of the resolutions of (\ref{skein2}), while the second one reduces the intersection number with $\Gamma$ and, hence, vanishes in the graded skein algebra.

\begin{figure}[h]
\centerline{\diag{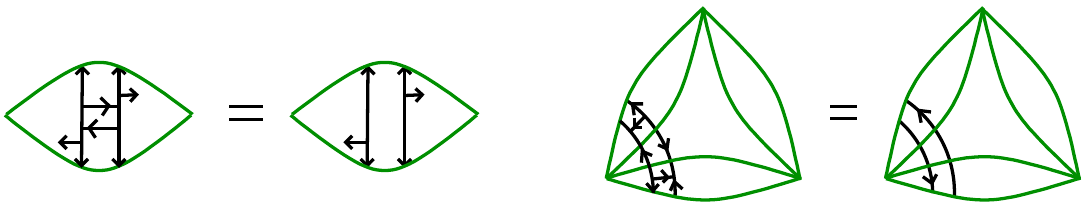}{4.2in}} 
\caption{$4$-gon eliminations. (Equalities in $\cal{GRS}(\hat F)$)}
\lb{f-4-gon-elimination}
\end{figure}

Although not essential for us here, one can easily prove that almost canonical webs have no $2$-gons (neither internal nor external), nor trivial components (contractible loops), nor external trigons.
In other words, almost canonical webs satisfy all conditions of being reduced, except possibly having  quadrigons. 

Let $W'$ be obtained now from $W_1$ by transposing two arcs in a face $T$ of the triangulation and by adding crossbars in the bigons on both sides of these arcs, as in Figure \ref{f-transposition}.

\begin{figure}[h]
\centerline{\diag{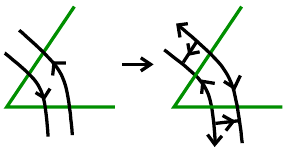}{1.5in}}
\caption{A transposition of arcs of opposite orientation.}
\lb{f-transposition}
\end{figure}

Then $W'$ is in an almost canonical position and equals to $W_1$ in $\cal{GRS}(\hat F).$ 
Let us repeat this process, until we obtain a web which coincides with $W_2$ in all faces of the triangulation. This web is in an almost canonical position and it equals to $W_1$ in $\cal{GRS}(\hat F).$ Let $W''$ be the web obtained by the elimination of all $4$-gons in the bigons. 
That elimination maintains its almost canonical position and does not affect its value in $\cal{GRS}(\hat F).$ Since $W''$ and $W_2$ coincide in all faces of the triangulation, the crossbar webs 
$W''\cap B$ and $W_2\cap B$ have coinciding signatures for internal every bigon $B$ of the triangulation. 
That is also the case for the external bigons since $W''$ and $W_2$ are trivial in these. 
Since the crossbar webs $W''\cap B$ and $W_2\cap B$ are minimal, they coincide by Lemma  \ref{l-red-signature}(3).
Consequently, $W''=W_2.$ 
\qed

%
\section{The coordinate monoid}
%

The coordinates of the canonical webs with respect to a triangulation $\Gamma$ form a subset $\cal C$ of $\Z_{\geq 0}^{\Gamma+\Gamma+
\cal F(\Gamma)},$
where $\cal F(\Gamma)$ denotes the set of faces of $\Gamma$.

The \underline{profile} of a canonical web $W$ is $\ve_W: \cal F(\Gamma)\to \{\pm 1\}$, where
$\ve_W(T)$ is the sign of the degree of the pyramid of $W$ in an ideal triangle $T$ of the ideal triangulation.
Denote the set of coordinates of all canonical webs of profile $\ve$ by $\cal S_\ve \subset \Z_{\geq 0}^{\Gamma+\Gamma+
{\cal F}(\Gamma)}$.

\bpro\lb{cor-Sve-fg}
$\cal S_\ve$ is a finitely generated additive submonoid of $\Z_{\geq 0}^{\Gamma+\Gamma+{\cal F}(\Gamma)}$
for every $\ve$.
\epro

\bpr 
Consider first the $7$ coordinates $(n_{01},n_{10},n_{12},n_{21},n_{20},n_{02},d)\in \Z_{\geq 0}^{7}$ for each face of $\Gamma$ for webs in $S_\ve,$ where $\ve(T)\cdot d$ is the degree of the pyramid in $T.$ Combined over all faces of $\Gamma$, they form an alternative set of web coordinates with values in $\Z_{\geq 0}^{7|{\cal F}(\Gamma)|}$.

If $T$ is bounded by $\gamma_{i_0}, \gamma_{i_1},\gamma_{i_2},$ oriented according to the orientation of $T$ then
these coordinates determine the intersection coordinates by (\ref{e-n2e}):
\[
\left(\begin{matrix}
e_{+,\gamma_{i_0}(W)}\cr e_{-,\gamma_{i_0}(W)}\cr e_{+,\gamma_{i_1}(W)}\cr e_{-,\gamma_{i_1}(W)}\cr e_{+,\gamma_{i_2}(W)}\cr e_{-,\gamma_{i_2}(W)}\cr r 
\end{matrix}\right)=\left(\begin{matrix}
0 & 1 & 0 & 0 & 1 & 0 & 1\\
1 & 0 & 0 & 0 & 0 & 1 & 0\\
1 & 0 & 0 & 1 & 0 & 0 & 1\\
0 & 1 & 1 & 0 & 0 & 0& 0\\
0 & 0 & 1 & 0 & 0 & 1& 1\\
0 & 0 & 0 & 1 & 1 & 0& 0\\
1 & -1 & 1 & -1 & 1 & -1& 0\\
\end{matrix}\right)\cdot \left(\begin{matrix} n_{01}\\ n_{10}\\ n_{12}\\ n_{21}\\n_{20}\\ n_{02}\\ d\end{matrix}\right),
\]
assuming the degree of the pyramid in $T$ is non-negative.
For the negative degree one needs to interchange ones and zeros in the first six entries of the last column.
The above numbers 
in $\Z_{\geq 0}^{7|{\cal F}(\Gamma)|}$ are realized by a web iff
for every internal edge $\gamma$ of the triangulation, the numbers $e_{\pm,\gamma(W)}$ coming from the two adjacent ideal triangles coincide. 
Such $7|{\cal F}(\Gamma)|$-tuples form an additive submonoid of $\Z_{\geq 0}^{7|{\cal F}(\Gamma)|}$. This monoid is finitely generated by Gordan's Lemma. $\cal S_\ve$ is the image of this monoid under the linear map given by the matrix above and, hence, a finitely generated (additive) monoid as well.
\epr

The empty web, which has all its coordinates zero is the unit of $\cal S_\ve.$

%
\section{Proof of Theorem \ref{main-gr}}
\lb{s_main_proof}
%


Note that the degree $d>0$ pyramid in a triangle $T$ is obtained by I-resolutions of all crossings in the stack of triads, $ST_d$, in Figure \ref{f-horizontal-triads}.
\begin{figure}[h]
\begin{center}
\begin{picture}(2.5,1.2)
\put(0,.6){\diag{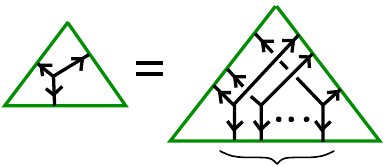}{2.7in}} 
\put(.5,.53){$d$}
\put(1.9,-.1){$d$}
\end{picture}
\end{center}
\caption{Stack of triads (sources), $SD_d$, in a horizontal formation.}
\lb{f-horizontal-triads}
\end{figure}
(Negative degree pyramid is obtained from stacking sink triads.)
Let $W$ be a non-elliptic web in $\hat F$ containing $T$ as a face of its triangulation. If $W\cap T= ST_d$ then
the geometric intersection number of $W$ with $\p T$ is $|ST_d\cap \p T|=3d.$
Note that on the other hand, any web $W'$ obtained by a vertexless resolution (as in Figure \ref{fig-cross-res}) of one of the crossings of $ST_d$ contains a vertex with two edges connecting it the same side of $T$. Hence, that web can be isotoped to reduce its intersection number with $\p T$.

\bcor
Let $W$ be a web in $\hat F$ with $ST_d$ in $T$ and $W'$ the same web with the pyramid $P_d$ in $T$ instead.
Then $W$ and $W'$ equal in $\cal{GRS}(\hat F)$ up to multiplication by a power of $q$.
\ecor

Let webs $W_1,W_2$ be reduced webs in canonical positions, of the same profile $\ve$.
Let us stack them up so that:
\bite
\item  in each face $T$ of the ideal triangulation:\\ 
\hspace*{.05in} (a) the corner arcs of $W_1$ and of $W_2$ are disjoint,\\
\hspace*{.05in} (b) the pyramids of $W_1$ and of $W_2$ intersect like in the Figure \ref{f-2honeycombs}.  (Hence, if they are of degrees $d_1$ and $d_2$ then they intersect $d_1\cdot d_2$ times.)
\item in each bigon, the vertical lines of $W_1$ and of $W_2$ are disjoint, the crossbars (horizontal lines) are disjoint and any vertical line of one intersects a horizontal line of the other at most once.
\eite

\begin{figure}[h]
\begin{center}
\begin{picture}(2,1.2)
\put(0,.6){\diag{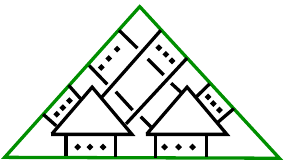}{2in}} 
\put(.56,.35){$P_{d_1}$}
\put(1.25,.35){$P_{d_2}$}
\end{picture}
\end{center}
\caption{Stacking pyramids}
\lb{f-2honeycombs}
\end{figure}

By (\ref{e-boundary-cross1}) and (\ref{e-boundary-cross2}), $W_1\cup W_2$ (as above) equals to $W_1\cdot W_2$
in $\cal{GRS}(\hat F)$ up to multiplication by a power of $q$.

Let $W_{12}$ be obtained by I-resolutions of all crossings of $W_1\cup W_2.$
Since $H_d$ can be obtained by I-resolutions of crossings in a horizontal stack of $d$ triads, 
the I-resolution of all vertices of the web in Fig. \ref{f-2honeycombs} yields a pyramid of degree $d_1+d_2,$
cf. Fig. \ref{f-honeycomb-addition}.

\begin{figure}[h]
\begin{center}
\begin{picture}(2,0.8)
\put(0,0.2){\diag{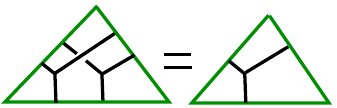}{2.4in}} 
\put(.45,.02){$d_1$}
\put(.8,.02){$d_2$}
\put(1.8,.01){$d_1$+$d_2$}
\end{picture}
\end{center}
\caption{Horizontal triad stacks.}
\lb{f-honeycomb-addition}
\end{figure}
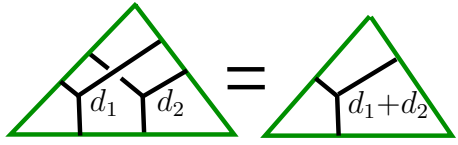

Hence, $W_{12}$ is non-elliptic and in canonical position in each triangle $T$ of the triangulation.

In principle, $W_{12}$ may contain however $4$-gons in bigons, cf. \ref{f-crossbar-webs-stacked}.
\begin{figure}[h]
\diag{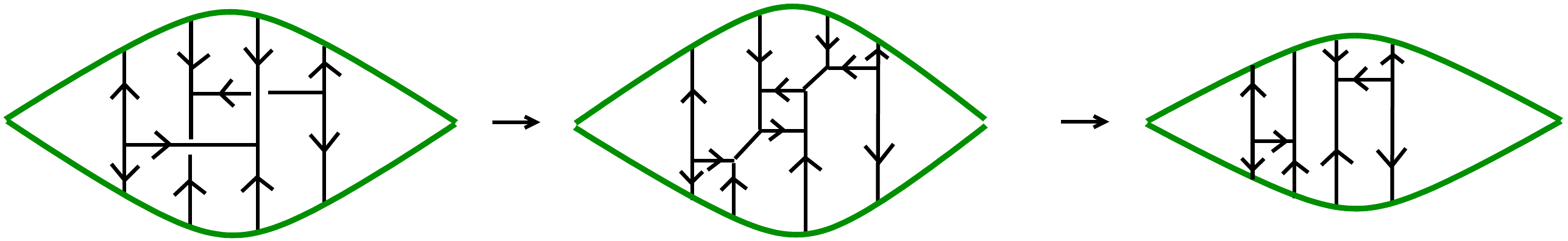}{4in} 
\caption{Left: Two crossbar webs stacked in a bigon. Middle: The result of I-resolving all crossings.
Right: The result of $4$-gon elimination.}
\lb{f-crossbar-webs-stacked}
\end{figure}

Each such $4$-gon in a bigon $B$ resolves into two horizontal arcs plus two vertical ones, cf. (\ref{skein2}).
However, the horizontal resolution leads to a web with a lower intersection number with $\p_{\pm} B.$
Let $\overline{W_{12}}$ be obtained from $W_{12}$ by resolving all $4$-gons in bigons vertically. Then by the above discussion,
\beq\lb{e-W1W2-12}
 W_1\cdot W_2 = q^{c}\cdot \overline{W_{12}}
\eeq
in $\cal{GRS}(\hat F)$ for some $c\in \frac{1}{3}\Z.$
Observe also that the coordinates of $\overline{W_{12}}$ are the sum of those of $W_1$ and $W_2$:
\beq\lb{e-W1W2-coord} 
(e_+(\overline{W_{12}}),e_-(\overline{W_{12}}),r(\overline{W_{12}}))=
\eeq
$$(e_+(W_1),e_-(W_1),r(W_1))+(e_+(W_2),e_-(W_2),r(W_2)\in \Z_{\geq 0}^{\Gamma+\Gamma+\cal F(\Gamma)}.$$

Now we are ready for:\vspace*{.1in}

\noindent{\bf Proof of Theorem \ref{main-gr}:}
By Proposition \ref{cor-Sve-fg}, $S_\ve$ is generated by a finite set $G_\ve\subset S_\ve$ for every $\ve: \cal F(\Gamma)\to \{\pm 1\}$.
Consider the union $G=\bigcup_\ve G_\ve$ of these generating sets for all $\ve$'s. We identify these elements with the reduced webs in $\hat F$ they correspond to through Theorem \ref{t-coordinates}. We claim that they generate ${\cal GRS}(\hat F)$. To this end we will show that every reduced web $W$ is a polynomial in those in $G$ by contradiction: Let $W$ be a web of smallest weight for which it is not the case. By Theorem \ref{t-main-canonical},  $W$ can be put into canonical position. Let $\ve$ be the profile of $W$. The coordinates of $W$ are either one of those in $G_\ve$, implying that $W\in G$, or they decompose 
$$(e_+(W),e_-(W),r(W))=(e_+(W_1),e_-(W_1),r(W_1))+(e_+(W_2),e_-(W_2),r(W_2)),$$
for some non-empty webs $W_1, W_2$, which by our assumption are polynomial expressions in elements of $G$. By (\ref{e-W1W2-12}) and (\ref{e-W1W2-coord}), $W$ equals $q^c\cdot W_1\cdot W_2$ in $\cal{GRS}(\hat F)$, for some $c\in \frac{1}{3}\Z$.
Since $W_1,W_2$ are of lower weight than $W$, that contradicts the assumption of $W$ being a lowest weight web which is not a polynomial expression in the webs in $G$.

%

\end{document}